# Character theory approach to Sato-Tate groups

Yih-Dar SHIEH

Abstract

In this article, we propose to use the character theory of compact Lie groups and their orthogonality relations for the study of Frobenius distribution and Sato-Tate groups. The results show the advantages of this new approach in several aspects. With samples of Frobenius ranging in size much smaller than the moment statistic approach, we obtain very good approximation to the expected values of these orthogonality relations, which give useful information about the underlying Sato-Tate groups and strong evidence of the correctness of the generalized Sato-Tate conjecture. In fact, $2^{10}$ to $2^{12}$ points provide satisfactory convergence. Even for $g = 2$, the classical approach using moment statistics requires about $2^{30}$ sample points to obtain such information.

## 1. *Introduction*

In [4], Fité, Kedlaya, Rotger and Sutherland study the limiting distributions of (the conjugacy classes of) the normalized Frobenius endomorphisms of abelian surfaces $A$ over number fields $K$, where the distributions are over primes of good reduction of $A/K$. Such distributions are expected to correspond to some closed subgroups $\mathrm{ST}_A$ of $\mathrm{USp}(4)$, on which the (conjugacy class of the) characteristic polynomial of a uniform random matrix gives the Frobenius distribution on $A$. The group $\mathrm{ST}_A$ is called the *Sato-Tate* group of $A$. They give a classification of Sato-Tate groups which shows that, up to conjugacy, there are exactly 52 groups which occur as Sato-Tate groups for suitable $A$ and $K$. They also exhibit examples of Jacobians of genus 2 hyperelliptic curves for each Sato-Tate group.

By the statistics of moments of the coefficients of the normalized characteristic polynomials of Frobenius of $A/K$, they empirically verified the expected Sato-Tate distribution. For each example curve, they compute sample points of Frobenius at primes $\mathfrak{p}$ of good reduction with norm $\|\mathfrak{p}\| \leq 2^{30}$. In the genus 2 case, these computations could be done practically by using the optimizations described in [10], which combines efficient point enumeration with generic group algorithms as discussed in Sutherland's PhD thesis [15], together with further improvements, especially an efficient implementation of the group operation in the Jacobian of curves, incorporated in the **smalljac** software library [16].

In [3], Fité and Sutherland study the Sato-Tate groups for the curves $y^2 = x^8 + c$ and $y^2 = x^7 - cx$. For these genus 3 curves $C$, they compute the Frobenius for primes $p \leq 2^{40}$ of good reduction of $C$, using efficient algorithms for curves in these families.

For generic hyperelliptic curves $C/\mathbb{Q}$ of genus $g$, Harvey's algorithm [5] computes the zeta function of the reduction $C_p$ of $C$ at $p$ for all primes $p \leq N$ of good reduction of $C$, where $N$ is a given bound. Its average complexity per prime is polynomial in $\log N$. Based on this work, in [6] and [7], Harvey and Sutherland present an efficient algorithm to compute the Hasse-Witt matrix of $C_p$ for all $p \leq N$, which gives the Frobenius characteristic polynomial $\chi_p$ modulo $p$. For $g \leq 3$, we can even determine $\chi_p$ by combining a generic group algorithm. This makes the computation up to $N = 2^{30}$ feasible for $g \leq 3$. However, in the study of Sato-Tate groups for $g = 3$, the results of moment statistics with $N = 2^{30}$ might not be satisfying.

---

In this article, instead of considering the moment statistics of the coefficients of the normalized Frobenius characteristic polynomials, we propose to use the orthogonality relations of the irreducible characters of the unitary symplectic group $\mathrm{USp}(2g)$ for the study of Sato-Tate groups in genus $g$. In Section 2, we first give an introduction to the question of Frobenius distributions. We then define the Sato-Tate group (Definition 1) and state the generalized Sato-Tate conjecture (Conjecture 1). This involves the notion of equidistribution (Definition 4), which is defined in Section 3, where we also recall some other notions from probability theory, and we present the orthogonality relations as expected values of certain random variables. In Section 4, we present a recursive algorithm (Algorithm 1) to compute the irreducible characters of $\mathrm{USp}(2g)$ based on the Brauer-Klimyk formula (Theorem 4.1), in terms of the coefficients of the normalized (real) Frobenius characteristic polynomial.

After introducing the background and necessary tools, we demonstrate the advantages of using orthogonality relations of irreducible characters through several examples in Section 5. In particular, in Example 2, we compare this new approach with the one using moment statistics. Example 5 gives a heuristic reason that our approach works very well whenever the orthogonality relations are given by small integers (e.g. the generic cases). We also propose a solution for non-generic cases, demonstrated in Example 3 and Example 4. A summary of these advantages are given in Section 6.

## 2. *Frobenius distributions and Sato-Tate groups*

In this section, we explain the two main objects that we study in this article: Frobenius distributions and Sato-Tate groups.

### 2.1. *Frobenius distribution*

Let $A/K$ be an abelian variety of dimension $g$, over a number field $K$. In almost all the examples in this article, $A = \mathrm{Jac}(C)$ is the Jacobian of some genus $g$ curve $C/K$ which has a $K$-rational point, and usually $K = \mathbb{Q}$.

Denote the set of all (finite) primes of $K$ by $M_K^0$. Let $S$ be the finite set of primes $\mathfrak{p}$ of bad reduction of $A$, and $O_{K,S}$ be the ring of $S$-integers of $K$. Let $\mathcal{A}$ be a model of $A$ over $O_{K,S}$, i.e. its special fiber $\mathcal{A} \times K$ is $A/K$. The set $\mathcal{P} = M_K^0 - S$ consists of primes of good reduction.

For each $\mathfrak{p} \in \mathcal{P}$, we obtain a reduction $\overline{A}_\mathfrak{p}$, which is an abelian variety over the residue field $k_\mathfrak{p} \simeq \mathbb{F}_q$ of $O_{K,S}$ at $\mathfrak{p}$, where $q := N_{K/\mathbb{Q}}(\mathfrak{p}) = \#k_\mathfrak{p}$. The characteristic polynomial $\widetilde{f}_\mathfrak{p}$ of the Frobenius action $\mathrm{Frob}_\mathfrak{p}$ on the rational Tate module $V_\ell(\overline{A}_\mathfrak{p}) = T_\ell(\overline{A}_\mathfrak{p}) \otimes_{\mathbb{Z}_\ell} \mathbb{Q}_\ell$ over $\mathbb{Q}_\ell$ is monic in $\mathbb{Z}[T]$ of degree $2g$. The Hasse-Weil theorem says that all the roots $\widetilde{\alpha}_1, \ldots, \widetilde{\alpha}_{2g}$ of $\widetilde{f}_\mathfrak{p}$ have absolute values $\sqrt{q}$, hence the normalized characteristic polynomial $f_\mathfrak{p} := \widetilde{f}_\mathfrak{p}\left(\sqrt{q}\,T\right)/q^g$, which has roots $\widetilde{\alpha}_1/\sqrt{q}, \ldots, \widetilde{\alpha}_{2g}/\sqrt{q}$, corresponds to a unique conjugacy class of the unitary symplectic group $\mathrm{USp}(2g)$. Fixing an embedding $\iota : \mathbb{Q}_\ell \hookrightarrow \mathbb{C}$, we have the normalized Frobenius action $\mathrm{NFrob}_\mathfrak{p} := \mathrm{Frob}_\mathfrak{p} \otimes \frac{1}{\sqrt{N(\mathfrak{p})}}$ on $V_\ell(\overline{A}_\mathfrak{p}) \otimes_{\mathbb{Q}_\ell} \mathbb{C}$, whose characteristic polynomial is $f_\mathfrak{p}$. This action is symplectic with respect to the Weil pairing, when we consider a fixed polarization of $A$. The Frobenius distribution (of $A/K$) is the distribution of the conjugacy class $[\mathrm{NFrob}_\mathfrak{p}]$ in $\mathrm{Cl}(\mathrm{USp}(2g))$ when $\mathfrak{p}$ varies over $S$, here $\mathrm{Cl}(G)$ is the set of conjugacy classes of a group $G$.

The Sato-Tate distributions and the (generalized) Sato-Tate conjecture (Conjecture 1) are concerned with the equidistribution of Frobenius. Consider the whole set $\{[\mathrm{NFrob}_\mathfrak{p}]\}_{\mathfrak{p}\in\mathcal{P}} \subseteq \mathrm{Cl}(\mathrm{USp}(2g))$. We search for a probability space $(G, \mathfrak{B}_G, \mu_G)^\dagger$, where $G \subset \mathrm{USp}(2g)$ is a compact Lie subgroup such that $\mathrm{Cl}(G)$ is the sample space in which $[\mathrm{NFrob}_\mathfrak{p}]$ live, and the probability measure $\mu_G$ is the (unique) normalized Haar measure on $G$, which is translation invariant. Serre propose a candidate $\mathrm{ST}_A$ of such group, which is called the Sato-Tate group of $A$, see Definition 1. The generalized Sato-Tate conjecture states that the distribution of $[\mathrm{NFrob}_\mathfrak{p}]$ is determined by the induced measure on $\mathrm{Cl}(G)$ from $\mu_G$. See Conjecture 1 for the precise statement.

$^\dagger$ Recall that the Borel $\sigma$-algebra $\mathfrak{B}_G$ on $G$ is the $\sigma$-algebra generated by the open subsets of $G$.

The definition of equidistribution (Definition 4) involves with the limits of sequences of sample statistics. In practice, it is impossible to gather information of $[\mathrm{NFrob}_\mathfrak{p}]$ for all $\mathfrak{p} \in \mathcal{P}$, we need to work with a sample, i.e. a chosen finite subset of $\{[\mathrm{NFrob}_\mathfrak{p}]\}_{\mathfrak{p}\in\mathcal{P}}$. In general, a sample is used to draw inferences and conclusions from itself to the whole set which we are concerned with. We usually compute $f_\mathfrak{p}$ for $\|\mathfrak{p}\| \le N$ for a chosen bound $N$. In the case $K = \mathbb{Q}$, we usually choose the first $n$ prime numbers in $\mathcal{P}$.

A conjugacy class $[\mathrm{NFrob}_\mathfrak{p}]$ is uniquely determined by its characteristic polynomial $f_\mathfrak{p}$. This way, the Frobenius distribution concerns the distribution of $f_\mathfrak{p}$ as $\mathfrak{p}$ varies over $\mathcal{P}$. In particular, we regard the map

$$\begin{aligned} \xi : \mathrm{Cl}(\mathrm{USp}(2g)) &\longrightarrow \mathbb{C}[T] \\ x &\longmapsto \mathrm{charpoly}(x) \end{aligned}$$

as a random variable[†]. For each $[\mathrm{NFrob}_\mathfrak{p}]$, we have $\xi([\mathrm{NFrob}_\mathfrak{p}]) = f_\mathfrak{p}$. We consider $\{f_\mathfrak{p}\}_{\mathfrak{p}\in\mathcal{P}}$, and the sample becomes the corresponding subset of $\{f_\mathfrak{p}\}_{\mathfrak{p}\in\mathcal{P}}$. We then study the distribution of $[\mathrm{NFrob}_\mathfrak{p}]$ via the sample statistics of the calculated sample $f_\mathfrak{p}$ .

From the functional equation $T^{2g}\widetilde{f}_\mathfrak{p}(1/T) = \widetilde{f}_\mathfrak{p}(qT)/q^g$ of the Weil polynomial $\widetilde{f}_\mathfrak{p}$, we obtain $f_\mathfrak{p}(T) = T^{2g} f_\mathfrak{p}(1/T)$, and the coefficients $a_i$ of $f_\mathfrak{p}$ satisfy $a_{2g-i} = a_i$. The normalized real Weil polynomial $g_\mathfrak{p} \in \mathbb{R}[T]$ is of degree $g$ satisfying $f_\mathfrak{p}(T) = T^g g_\mathfrak{p}(T + 1/T)$. Since the characteristic polynomials (of different types) of Frobenius are all of their own importance, we fix the following notations:

$$\widetilde{f}_\mathfrak{p}(T) = T^{2g} - \widetilde{a}_1 T^{2g-1} + \widetilde{a}_2 T^{2g-2} - \ldots + \widetilde{a}_2 q^{g-2} T^2 - \widetilde{a}_1 q^{g-1} T + q^g$$

and

$$f_\mathfrak{p}(T) = T^{2g} - a_1 T^{2g-1} + a_2 T^{2g-2} - \ldots + a_2 T^2 - a_1 T + 1,$$

where $a_i = \widetilde{a}_i / \sqrt{q}^{\,i}$. Letting $t_i = \alpha_i + \alpha_i^{-1}$, we obtain

$$f_\mathfrak{p}(T) = \prod_{i=1}^{g} (T - \alpha_i)(T - \alpha_i^{-1}) = \prod_{i=1}^{g} (T^2 - t_i T + 1).$$

Finally, we define

$$g_\mathfrak{p}(T) = \prod_{i=1}^{g} (T - t_i) = T^g - s_1 T^{g-1} + s_2 T^{g-2} - \ldots + (-1)^{g-1} s_{g-1} T + (-1)^g s_g,$$

where $s_i = \mathrm{sym}(t_1, \ldots, t_g)$, the $i$-th elementary symmetric function. For $x \in \mathrm{USp}(2g)$ or $\mathrm{Cl}(\mathrm{USp}(2g))$, we write $\widetilde{f}_x, f_x$ and $g_x$ for its characteristic polynomial, normalized characteristic polynomial and normalized real characteristic polynomial respectively. Instead of working with the random element $\xi$ above, which has values in $\mathbb{C}[T]$, we consider the random variables (for $1 \le i \le g$)

$$\begin{aligned} a_i : \mathrm{Cl}(\mathrm{USp}(2g)) &\longrightarrow \mathbb{C} \\ x &\longmapsto (-1)^i \times \text{ the coefficient of } T^{2g-i} \text{ in } f_x \end{aligned} \tag{2.1}$$

[†] See Section 3.

The authors in [4] use the sample moment statistics of $a_i$ to study the Frobenius distributions. Instead of $a_i$, one may use the random variables (for $1 \leq i \leq g$)

$$\begin{aligned} s_i : \mathrm{Cl}(\mathrm{USp}(2g)) &\longrightarrow \mathbb{C} \\ x &\longmapsto (-1)^i \times \text{ the coefficient of } T^{g-i} \text{ in } g_x \end{aligned} \tag{2.2}$$

since $g_x$ determines $f_x$ and vice versa. However, we will use the orthogonality relations of the irreducible characters of $\mathrm{USp}(2g)$ in Section 5 to study the Frobenius distribution and Sate-Tate groups, and we demonstrate the advantages of this new approach.

Remark 1. Let $s_0 = a_0 = 1$. It is easy to prove that $a_j = \sum_{i=0}^{g} c_{i,j} s_i$, where $c_{i,j} \in \mathbb{Z}$ (depending on $g$) is determined by the following recurrence relation (for all $j \in \mathbb{N}$):

$$\begin{aligned} c_{g,j} &= 0, \text{ if } j \neq i \\ c_{g,g} &= 1 \\ c_{i,j} &= c_{i+1,j-1} + c_{i+1,j+1} \end{aligned}$$

We have a closed formula (for $0 \leq i, j \leq g$)

$$c_{i,j} = \frac{1 + (-1)^{i+j}}{2} \binom{g-i}{g - \frac{i+j}{2}}.$$

The expressions $s_j = \sum_{i=0}^{g} d_{i,j} a_i$ of $s_j$ in $a_i$ is given by

$$d_{i,j} = \left(\mathbf{i}^{i-j} + \mathbf{i}^{j-i}\right) \frac{g-i}{2g-i-j} \binom{g - \frac{i+j}{2}}{\frac{j-i}{2}},$$

where $d_{g,g} = 1$ and $\mathbf{i} \in \mathbb{C}$ is the imaginary unit.

## 2.2. *Sato-Tate groups*

We refer to Serre's book [11], a lecture note [9] of Kedlaya or of Sutherland [14] for the definition of the Sato-Tate group (Definition 1) and the generalized Sato-Tate conjecture (Conjecture 1).

Let $A/K$ be as in Section 2.1. We fix a prime number $\ell$, and the set $S$ of primes $\mathfrak{p}$ is as in Section 2.1, but also including those $\mathfrak{p}$ lying over $\ell$, which is again a finite set. Let $\mathcal{P} = M_K^0 - S$.

Let $K_{A,\ell} = K\left(A\left[\ell^\infty\right]\right) \subseteq \overline{\mathbb{Q}}$ be the $\ell^\infty$-division field of $A$ and $\rho_{A,\ell} : \mathrm{Gal}(\overline{K}/K) \to \mathrm{Aut}(T_\ell(A))$ be the $\ell$-adic representation attached to the abelian variety $A/K$. It factors through the quotient:

$$\begin{array}{ccc} \mathrm{Gal}(\overline{K}/K) & \longrightarrow & \mathrm{Gal}(K_{A,\ell}/K) \\ & {\scriptstyle \rho_{A,\ell}} \searrow \quad \swarrow {\scriptstyle \rho_{A,\ell}} & \\ & \mathrm{Aut}(T_\ell(A)) & \end{array}$$

Consider†

$$\begin{aligned} \mathrm{Frob} : \mathcal{P} &\longrightarrow \mathrm{Gal}(K_{A,\ell}/K) \xrightarrow{\rho_{A,\ell}} \mathrm{Aut}(T_\ell(A) \otimes_{\mathbb{Z}_\ell} \mathbb{Q}_\ell) \\ \mathfrak{p} &\longmapsto \sigma_{\mathfrak{P}} \longmapsto F_{\mathfrak{P}} \end{aligned}$$

† Recall that those primes $\mathfrak{p}$ dividing $\ell$ are excluded from the set $\mathcal{P}$, as mentioned in the beginning of Section 2.2.

where $\sigma_{\mathfrak{P}}$ is the Frobenius element of a choice of place $\mathfrak{P}$ over $\mathfrak{p}$[‡], and $F_{\mathfrak{P}}$ is induced from the action of $\sigma_{\mathfrak{P}}$ on $A(K_{A,\ell})$. Subject to the choices of places $\mathfrak{P}$ over primes $\mathfrak{p}$, the map is well-defined because $\mathfrak{p} \nmid \ell$ is unramified[§] in $K_{A,\ell}$. Different choices of $\mathfrak{P}$ determine conjugate Frobenius elements $\sigma_{\mathfrak{P}}$, and hence conjugate actions $\mathrm{Frob}(\mathfrak{p}) = F_{\mathfrak{P}}$.

We have canonical isomorphisms (induced from the reduction modulo $\mathfrak{p}$) such that the following diagram is commutative:

$$\begin{array}{ccc}
T_\ell(A) \otimes_{\mathbb{Z}_\ell} \mathbb{Q}_\ell & \xrightarrow{\ \sim\ } & T_\ell(\bar{A}_{\mathfrak{p}}) \otimes_{\mathbb{Z}_\ell} \mathbb{Q}_\ell \\
{\scriptstyle \mathrm{Frob}(\mathfrak{p})}\Big\uparrow & & \Big\uparrow{\scriptstyle \mathrm{Frob}_{\mathfrak{p}}} \\
T_\ell(A) \otimes_{\mathbb{Z}_\ell} \mathbb{Q}_\ell & \xrightarrow{\ \sim\ } & T_\ell(\bar{A}_{\mathfrak{p}}) \otimes_{\mathbb{Z}_\ell} \mathbb{Q}_\ell
\end{array}$$

The Frobenius actions $\mathrm{Frob}_{\mathfrak{p}}$ on different spaces $V_\ell(\bar{A}_{\mathfrak{p}})$ are then realized by the actions $\mathrm{Frob}(\mathfrak{p})$ on a common space $V_\ell(A) = T_\ell(A) \otimes_{\mathbb{Z}_\ell} \mathbb{Q}_\ell$.

Fix a polarization for $A/K$, which gives the Weil pairing[†] $e_\ell$ on the rational $\ell$-adic Tate module $V_\ell(A)$, making $V_\ell(A)$ a symplectic vector space over $\mathbb{Q}_\ell$. The representation $\rho_{A,\ell}$ is symplectic, i.e. $e_\ell(\rho_{A,\ell}(\sigma) \cdot v, \rho_{A,\ell}(\sigma) \cdot w) = e_\ell(v, w)^\sigma$ for all $\sigma \in \mathrm{Gal}(\overline{K}/K)$ and $(v, w) \in V_\ell(A) \times V_\ell(A)$. By fixing a $e_\ell$-symplectic basis for $V_\ell(A)$, we obtain $\rho_{A,\ell} : \mathrm{Gal}(\overline{K}/K) \to \mathrm{GSp}(2g, \mathbb{Q}_\ell)$. Let $G_\ell \subseteq \mathrm{GSp}(2g, \mathbb{Q}_\ell)$ be the Zariski closure of $\rho_{A,\ell}(\mathrm{Gal}(\overline{K}/K))$ and $G_\ell^1 = G_\ell \cap \mathrm{Sp}(2g, \mathbb{Q}_\ell)$.

Definition 1. Choose an embedding $\iota : \mathbb{Q}_\ell \hookrightarrow \mathbb{C}$. Let $G^1 = G_\ell^1 \otimes_\iota \mathbb{C} \subseteq \mathrm{Sp}(2g, \mathbb{C})$. The *Sato-Tate group* $\mathrm{ST}_A$ of $A$ is a maximal compact Lie subgroup of $G^1$ contained in $\mathrm{USp}(2g)$.

We are concerned only with the conjugacy class $[\mathrm{NFrob}_{\mathfrak{p}}]$ in $\mathrm{Cl}(\mathrm{USp}(2g))$, which is given by $[\mathrm{NFrob}(\mathfrak{p})]$ where $\mathrm{NFrob}(\mathfrak{p}) := \mathrm{Frob}(\mathfrak{p}) \otimes \frac{1}{\sqrt{N(\mathfrak{p})}}$ on $V_\ell(A) \otimes_\iota \mathbb{C}$. It is generally expected that $\mathrm{NFrob}(\mathfrak{p})$ is conjugate to an element in $\mathrm{ST}_A$, hence $[\mathrm{NFrob}(\mathfrak{p})]$ is in the image $\mathrm{Cl}_A := \mathrm{Cl}(\mathrm{ST}_A) \to \mathrm{Cl}(\mathrm{USp}(2g))$.

Conjecture 1 (Generalized Sato-Tate Conjecture). For each positive integer $N$, let $S_N = \{[\mathrm{NFrob}(\mathfrak{p})]\}_{\|\mathfrak{p}\| \le N}$, which is a finite subset in $\mathrm{Cl}_A$. Then the sequence $(S_N)_{N=1}^\infty$ is equidistributed with respect to the induced measure of the Haar measure of $\mathrm{ST}_A$ on $\mathrm{Cl}_A$. See Definition 4 for the definition of equidistribution.

## 3. *Random variables, moments and equidistribution*

In this section, we recall some notions from probability theory. Let $(X, \Sigma, \mu)$ be a probability space and $\xi : X \to \mathbb{C}$ be a random variable, i.e. a measurable function on $X$ with respect to the $\sigma$-algebra $\Sigma$ on $X$ and the usual Lebesgue measure on $\mathbb{C}$.

---

‡ For the ramification theory for infinite Galois extensions, see [17, pp. 332–336, Appendix, §2]. For more details, see [8, §6]

§ This is a consequence of Néron-Ogg-Shafarevich criterion. See [12, Thm. 1].

† Here we fix a polarization on $A$ such that the corresponding Weil pairing $e_\ell$ is non-degenerate and skew symmetric.

Definition 2. The *expectation* $E[\xi]$ of the random variable $\xi$ is

$$E[\xi] = \int_{x\in X} \xi(x)\mu(dx)\cdot$$

The *n-th moment* $M_n[\xi]$ of $\xi$ is the expectation of $\xi^n$.

Definition 3. Given a sample $S \subseteq X$, i.e. a finite subset of $X$, the *n-th sample moment* of $\xi$ for the sample $S$ is

$$M_{n,S}[\xi] = \frac{1}{|S|}\sum_{x\in S}\xi^n(x)\cdot$$

The sample moment statistics $M_{n,S}$ are used to provide information of the probability distribution $(X,\mu)$, when it is unknown at the first place, or to give the empirical evidence for a conjectural distribution. In general, for a random sample $S$ whose size is sufficiently large, we expect $M_{n,S}[\xi]$ to be a good estimation of $M_n[\xi]$. The notion of equidistribution is based on this idea.

Let $(X,\mu)$ be a probability space where $X$ is a metric space and we use the Borel $\sigma$-algebra on $X$. Every continuous function $\xi$ is a measurable function on $X$, and thus a random variable.

Definition 4 (Equidistribution). Let $I$ be a totally ordered set (usually, $\mathbb{N}$ or $\mathbb{R}$). Let $(S_i)_{i\in I}$ be a family of finite subsets of $(X,\mu)$ satisfying $S_i \subseteq S_j$ if $i \le j$. The family $(S_i)$ is said to be *equidistributed with respect to the probability measure* $\mu$ if the following condition holds: For any bounded continuous function $\xi : X \to \mathbb{C}$, we have

$$\lim_{i\to\infty}\frac{1}{|S_i|}\sum_{x\in S_i}\xi(x) = \int_X \xi(x)\mu(dx).$$

A sequence $(x_k)_{k=1}^{\infty}$ in $(X,\mu)$ is said to be equidistributed with respect to $\mu$ if the family $\left(S_i = \{x_k\}_{k=1}^{i}\right)$ is equidistributed.

The equidistribution means that the sample mean $M_{1,S_i}[\xi]$ of the sample $S_i$ converges to the expected value $E[\xi] = M_1[\xi]$ of the random variable $\xi$. In particular, we have $M_{n,S_i}[\xi] \to M_n[\xi]$ for all higher moments of $\xi$.

Let $G$ be a compact Lie group and $\mu_G$ be its (normalized) Haar measure, which makes $(G,\mu_G)$ a probability space. Each (virtual) character $\chi : G \to \mathbb{C}$ can be regarded as a random variable. For $G \subseteq \mathrm{USp}(2g)$, one can consider the restrictions of $a_i$, $s_i$ and $\chi_i$, or more generally $\chi$, on $G$, where $a_i$ and $s_i$ are defined in Equation (2.1) and Equation (2.2), and $\chi_i$ (or $\chi$) are the (fundamental$^\dagger$) irreducible characters of $\mathrm{USp}(2g)$. The sample moment statistics of these random variables over a sample of Frobenius are used to obtain the conjectural Sato-Tate group, or to verify empirically the generalized Sato-Tate conjecture.

In this article, instead of using these moment statistics, we propose to use the orthogonality relations on $G$ of the irreducible characters of $\mathrm{USp}(2g)$. For two irreducible characters $\chi_\lambda$ and $\chi_\nu$ of $\mathrm{USp}(2g)$, we consider

$$\langle\chi_\lambda,\chi_\nu\rangle_G = E\left[\chi_\lambda\overline{\chi_\nu}|_G\right] = \int_{x\in G}\chi_\lambda\overline{\chi_\nu}(x)\,\mu_G(dx).$$

$^\dagger$These are the irreducible characters $\mathrm{USp}(2g)$ corresponding to "the" fundamental dominant weights $\varpi_i$, defined in Section 4, which are themselves determined by a choice of simple positive roots of $\mathrm{USp}(2g)$.

If $G = \mathrm{USp}(2g)$, we have $\langle \chi_\lambda, \chi_\nu \rangle = \delta_{\lambda,\nu}$ from Schur orthogonality, where $\delta_{\lambda,\nu}$, is the Kronecker delta symbol. In general, one needs to consider the branching rules from $\mathrm{USp}(2g)$ to $G$ to obtain the expected value of the random variable $\xi = \chi_\lambda \overline{\chi_\nu}|_G$ on $G$.

## 4. *Character Theory of* USp(2g)

In this section, we present a recursive method to compute the irreducible characters of $\mathrm{USp}(2g)$ based on the Brauer-Klimyk formula. We give a minimal background to introduce the notation and to present the results, and we refer to [1], in particular, Chapter 18–22, for the general theory of compact Lie groups, and their irreducible representations and characters.

We fix an embedding of the *unitary symplectic group*:

$$\mathrm{USp}(2g) = \{x \in \mathrm{GL}(2g, \mathbb{C}) \,|\, x^t J x = J \text{ and } \overline{x}^t x = I_{2g}\}, \ J = \begin{bmatrix} 0 & I_g \\ -I_g & 0 \end{bmatrix}.$$

We choose a maximal torus $T$ for $\mathrm{USp}(2g)$, which is given by diagonal matrices of the form:

$$u = \begin{bmatrix} u_1 & & & & & \\ & \ddots & & & & \\ & & u_g & & & \\ & & & \overline{u}_1 & & \\ & & & & \ddots & \\ & & & & & \overline{u}_g \end{bmatrix}, \ u_i \in \mathrm{U}(1).$$

A *weight* is a continuous homomorphism $\lambda : T \to \mathbb{C}^\times$. For $1 \leq i \neq j \leq g$, let $\alpha_{i,j} : T \to \mathbb{C}^\times$, $u \mapsto \frac{u_i}{u_j}$. For $1 \leq k \leq g$, we also define $\alpha_{k,k} : T \to \mathbb{C}^\times$, $u \mapsto u_k^2$. Our choice of a set of simple positive roots is $\{\alpha_k\}_{1 \leq k \leq g}$, where $\alpha_k = \alpha_{k,k+1}$ if $k \leq g-1$ and $\alpha_g = \alpha_{g,g}$. Under this choice, the fundamental dominant weights are $\varpi_k : T \to \mathbb{C}^\times$, $u \mapsto \prod_{i=1}^k u_i$, which form a basis $\varpi$ of the weight lattice $\Lambda$. Each dominant weight $\lambda$ is of the form $\sum_{i=1}^g n_i \varpi_i$ with all $n_i \in \mathbb{N}$. We work with the coordinate $[\lambda]_\varpi = (n_1, \ldots, n_g)$, and we define the *unweighted degree* of the dominant weight $\lambda$ to be $\sigma(\lambda) = \sum_{i=1}^g n_i$.

For a compact connected semisimple Lie group $G$, the theorem of the highest weight tells us that there is a one-to-one correspondence between the dominant weights of $G$ and the finite dimensional irreducible representations of $G$, and the irreducible character $\chi_\lambda$ of the representation $\rho_\lambda$ for a dominant weight $\lambda$ is given by the Weyl character formula. However, this is not suitable for efficient computation. One of the reasons is that it involves the explicit action for each element in the Weyl group, which is usually a huge group. Furthermore, unlike the recursive Algorithm 1 based on the Brauer-Klimyk formula, it loses the advantage of using the previously computed results, when the computation of (a sequence of) irreducible characters is concerned, rather than a single one. See the author's thesis [13, p. 96, Prop 4.68] for a discussion on the average time complexity per character of Algorithm 1.

Let $\Lambda_+$ be the set of dominant weights, $\rho = \sum_{i=1}^g \varpi_i$ be the Weyl vector, and $W = N(T)/T$ be the Weyl group of $G$.

THEOREM 4.1 (Brauer, Klimyk. See [1, p. 185, Prop. 22.9]). *Let $\lambda \in \Lambda_+$ and $\nu \in \Lambda$. There is $w \in W$ such that $\eta_\nu = w(\lambda + \rho + \nu) \in \Lambda_+$. The point $\eta_\nu$ is uniquely determined. If $\eta_\nu$ is on the boundary of $\Lambda_+$, we define $\xi_\nu = 0$. Otherwise $w$ is also uniquely determined, $\eta_\nu - \rho \in \Lambda_+$, and we define $\xi_\nu = \det(w)\chi_{\eta_\nu - \rho}$. For a dominant weight $\mu$ for which the weight decomposition is $\chi_\mu|_T = \sum m(\nu)\nu$, we have*

$$\chi_\mu \chi_\lambda = \sum_\nu m(\nu)\xi_\nu.$$

The Brauer-Klimyk formula can be turned into a recursive algorithm to compute the irreducible characters. For $G = \mathrm{USp}(2g)$, it is done in the author's thesis [13], where we devote some effort to prove the termination of the algorithm. We present the algorithm itself in Algorithm 1. Here, $\chi_i$ are the fundamental irreducible characters corresponding to $\varpi_i$ and $\chi_0 = 1$. For $x \in \mathrm{USp}(2g)$, one can consider its characteristic polynomial $f_x$ and its real characteristic polynomial $g_x$ as in Section 2.1, and their coefficients $a_i$ and $s_i$ (as functions in $x$). The relations between $\chi_i$, $s_i$ and $a_i$ are given in Lemma 4.2 and Corollary 4.3. These results should be well-known, however a proof is given in [13, §4.5].

**Algorithm 1** Compute Irreducible Characters of USp$(2g)$ in $\mathbb{Z}[\chi_1, \ldots, \chi_g]$

**def** CHI($x$) — # $\chi_\lambda$ for $[\lambda]_\varpi = x \in \mathbb{N}^g$
  **if** $x \notin \mathbb{N}^g$ **:** — # $\lambda$ should be dominant
    **return** 0
  **if** $\sigma(x) = 0$ **:** — # $x = (0, \ldots, 0)$
    **return** 1
  Find $1 \leq l \leq g$ such that $x_l \geq 1$
  **if** $\sigma(x) = 1$ **:** — # $x = e_l$
    **return** the symbol $\chi_l$ — # Recursive computing
  Set $\tilde{\chi} = \sum_{\nu \neq \varpi_l} \det(w)\, m(\nu)\, \mathrm{CHI}([w((\lambda - \varpi_l) + \rho + \nu) - \rho]_\varpi)$
  **return** CHI$(x - e_l)$ CHI$(e_l)$ $-\tilde{\chi}$

LEMMA 4.2. *We have $\chi_j = \sum_{i=0}^{g} c_{i,j} s_i$, where $c_{i,j} \in \mathbb{Z}$ (depending on $g$) is determined by the following recurrence relation:*

$$\begin{aligned} c_{g,j} &= 0, \text{ if } j \neq i \\ c_{g,g} &= 1 \\ c_{i,j} &= 0, \text{ if } j > g \\ c_{i,j} &= c_{i+1,j-1} + c_{i+1,j+1} \end{aligned}$$

*We have a closed formula*

$$c_{i,j} = \frac{1 + (-1)^{i+j}}{2} \left( \frac{2(g+1-j)}{2(g+1) - (i+j)} \right) \binom{g-i}{g - \frac{i+j}{2}}.$$

*The expressions $s_j = \sum_{i=0}^{g} d_{i,j} \chi_i$ of $s_j$ in $\chi_i$ is given by*

$$d_{i,j} = \frac{1}{2} \left( \mathbf{i}^{i-j} + \mathbf{i}^{j-i} \right) \binom{g - \frac{i+j}{2}}{\frac{j-i}{2}}.$$

COROLLARY 4.3. *We have $\chi_0 = a_0$, $\chi_1 = a_1$, and $\chi_i = a_i - a_{i-2}$ for $2 \leq i \leq g$.*

EXAMPLE 1. For $g = 2$ and $g = 3$, we have

| $\lambda$ | $\chi_\lambda$ in terms of | | |
|---|---|---|---|
| | $\chi_i$ | $s_i$ | $a_i$ |
| $(0,0)$ | $\chi_0$ | $s_0$ | $a_0$ |
| $(1,0)$ | $\chi_1$ | $s_1$ | $a_1$ |
| $(0,1)$ | $\chi_2$ | $s_2+1$ | $a_2-1$ |
| $(2,0)$ | $\chi_1^2-\chi_2-1$ | $s_1^2-s_2-2$ | $a_1^2-a_2$ |
| $(1,1)$ | $\chi_1\chi_2-\chi_1$ | $s_1s_2$ | $a_1a_2-2a_1$ |
| $(0,2)$ | $\chi_2^2-\chi_1^2+\chi_2$ | $s_2^2-s_1^2+3s_2+2$ | $a_2^2-a_1^2-a_2$ |
| $(3,0)$ | $\chi_1^3-2\chi_1\chi_2-\chi_1$ | $s_1^3-2s_1s_2-3s_1$ | $a_1^3-2a_1a_2+a_1$ |
| $(2,1)$ | $\chi_1^2\chi_2-\chi_2^2-\chi_1^2-\chi_2+1$ | $s_1^2s_2-s_2^2-3s_2-1$ | $a_1^2a_2-a_2^2-2a_1^2+a_2+1$ |
| $(1,2)$ | $\chi_1\chi_2^2-\chi_1^3+\chi_1$ | $s_1s_2^2-s_1^3+2s_1s_2+2s_1$ | $a_1a_2^2-a_1^3-2a_1a_2+2a_1$ |
| $(0,3)$ | $\chi_2^3-2\chi_1^2\chi_2+2\chi_2^2+\chi_1^2-1$ | $s_2^3-2s_1^2s_2+5s_2^2-s_1^2+7s_2+2$ | $a_2^3-2a_1^2a_2-a_2^2+3a_1^2-a_2$ |

$g=2$

| $\lambda$ | $\chi_\lambda$ in terms of | | |
|---|---|---|---|
| | $\chi_i$ | $s_i$ | $a_i$ |
| $(0,0,0)$ | $\chi_0$ | $s_0$ | $a_0$ |
| $(1,0,0)$ | $\chi_1$ | $s_1$ | $a_1$ |
| $(0,1,0)$ | $\chi_2$ | $s_2+2$ | $a_2-1$ |
| $(0,0,1)$ | $\chi_3$ | $s_3+s_1$ | $a_3-a_1$ |
| $(2,0,0)$ | $\chi_1^2-\chi_2-1$ | $s_1^2-s_2-3$ | $a_1^2-a_2$ |
| $(1,1,0)$ | $\chi_1\chi_2-\chi_3-\chi_1$ | $s_1s_2-s_3$ | $a_1a_2-a_1-a_3$ |
| $(1,0,1)$ | $\chi_1\chi_3-\chi_2$ | $s_1s_3+s_1^2-s_2-2$ | $a_1a_3-a_1^2-a_2+1$ |
| $(0,2,0)$ | $\chi_2^2-\chi_1\chi_3-\chi_1^2+\chi_2$ | $s_2^2-s_1s_3-2s_1^2+5s_2+6$ | $a_2^2-a_1a_3-a_2$ |
| $(0,1,1)$ | $\chi_2\chi_3-\chi_1\chi_2+\chi_3$ | $s_2s_3+3s_3+s_1$ | $a_2a_3-2a_1a_2+a_1$ |
| $(0,0,2)$ | $\chi_3^2-\chi_2^2+\chi_1\chi_3$ | $s_3^2-s_2^2+3s_1s_3+2s_1^2-4s_2-4$ | $a_3^2-a_2^2-a_1a_3+2a_2-1$ |

$g=3$

## 5. *Explicit computations*

### 5.1. *General framework*

In our study of Sato-Tate groups, the following objects are given or computed before the computation of Frobenius characteristic polynomials:

- An abelian variety $A$ (or a curve $C$) of dimension $g$ (or of genus $g$) over a number field $K$.
- A compact connected Lie subgroup $G$ of $\mathrm{USp}(2g)$ which we know contains the Sato-Tate group of $A/K$. Usually, $G=\mathrm{USp}(2g)$.
- A conjectural Sato-Tate group $H\subseteq G$ for $A/K$.
- A finite subset $S$ of the set $\mathcal{P}$ of primes $\mathfrak{p}$ of good reduction of $A/K$.
- A finite subset $I$ of the dominant weights of $G$. For $G=\mathrm{USp}(2g)$, we usually choose $I=I_d=\{\lambda\in\Lambda_+\mid\sigma(\lambda)\leq d\}$ for some positive integer $d$.

For each $\mathfrak{p}\in S$, we compute the normalized real Weil polynomial of $A/K$ at $\mathfrak{p}$, recorded by its coefficients $F_\mathfrak{p}=(s_1,s_2,\ldots,s_g)_\mathfrak{p}$. For any two dominant weights $\lambda$ and $\mu$ of $G$ which are in $I$, we compute the sample mean $M_{1,S}\left[\chi_\lambda\chi_\mu|_H\right]$ of the random variable $\chi_\lambda\chi_\mu$ on $H$.

In Example 5, where we study the heuristic behavior in the sample size, and in other examples in [13] regarding the heuristic behaviors in the genus $g$ and in the number of irreducible characters used, we denote the difference between the sample mean and the expected value of

$\chi_\lambda \chi_\mu$ [†] by

$$\begin{aligned} \operatorname{err}(H,S,\lambda,\mu) &= M_{1,S}\left[\chi_\lambda\chi_\mu|_H\right] - E\left[\chi_\lambda\chi_\mu|_H\right] \\ &= \frac{1}{|S|}\sum_{\mathfrak{p}\in S}\chi_\lambda(F_\mathfrak{p})\chi_\mu(F_\mathfrak{p}) - \langle\chi_\lambda,\chi_\mu\rangle_H. \end{aligned} \tag{5.1}$$

Finally, we compute $\operatorname{Err}(H,S,I) = \max_{\lambda,\mu\in I} \operatorname{err}(H,S,\lambda,\mu)$.

### 5.2. *Examples*

Example 2. In this example, we compare the moment statistics approach with our new approach using the orthogonality relations of irreducible characters. We take $C : y^2 = x^5 + x + 1$. Its conjectural Sato-Tate group is $H = G = \mathrm{USp}(4)$. In this example, we consider the first $N$ primes $p$ of good reduction of $C$, rather than those $p \leq N$. The moments of $a_1$ and $a_2$ are given in the columns with $N = \infty$. For $n \geq 5$, even with $2^{16}$ sample points, we don't obtain useful approximations of $M_n[a_1]$ and $M_n[a_2]$.

| $n$ | $N = 2^{12}$ | $N = 2^{16}$ | $N = \infty$ |
|---|---|---|---|
| 1 | 0.002 | 0.006 | 0 |
| 2 | 0.984 | 0.996 | 1 |
| 3 | 0.046 | −0.001 | 0 |
| 4 | 2.833 | 2.970 | 3 |
| 5 | 0.196 | −0.128 | 0 |
| 6 | 12.306 | 13.743 | 14 |
| 7 | 0.397 | −1.487 | 0 |
| 8 | 66.441 | 81.446 | 84 |
| 9 | −3.853 | −14.304 | 0 |
| 10 | 409.298 | 565.972 | 594 |

$M_n[a_1]$

| $n$ | $N = 2^{12}$ | $N = 2^{16}$ | $N = \infty$ |
|---|---|---|---|
| 1 | 0.989 | 0.999 | 1 |
| 2 | 1.964 | 1.992 | 2 |
| 3 | 3.815 | 3.966 | 4 |
| 4 | 9.250 | 9.853 | 10 |
| 5 | 23.747 | 26.423 | 27 |
| 6 | 67.907 | 79.611 | 82 |
| 7 | 205.367 | 257.730 | 268 |
| 8 | 658.293 | 893.546 | 940 |
| 9 | 2192.789 | 3257.407 | 3476 |
| 10 | 7550.758 | 12387.749 | 13448 |

$M_n[a_2]$

Table 3: $C : y^2 = x^5 + x + 1$, $H = G = \mathrm{USp}(4)$

Now we use the orthogonality relations of the irreducible characters of USp(4). We take the first six irreducible characters for $g = 2$ in Example 1, and we denote them by $\chi_i$ for $0 \leq i \leq 5$. We expect to see orthonormal relations between these $\chi_i$.

| | $\chi_0$ | $\chi_1$ | $\chi_2$ | $\chi_3$ | $\chi_4$ | $\chi_5$ |
|---|---|---|---|---|---|---|
| $\chi_0$ | 1.000 | −0.037 | 0.003 | 0.004 | −0.021 | −0.050 |
| $\chi_1$ | −0.037 | 1.007 | −0.058 | −0.095 | −0.057 | −0.017 |
| $\chi_2$ | 0.003 | −0.058 | 0.954 | −0.006 | −0.091 | −0.038 |
| $\chi_3$ | 0.004 | −0.095 | −0.006 | 0.928 | −0.054 | −0.071 |
| $\chi_4$ | −0.021 | −0.057 | −0.091 | −0.054 | 0.879 | −0.075 |
| $\chi_5$ | −0.050 | −0.017 | −0.038 | −0.071 | −0.075 | 0.947 |

Orthogonality relations with $N = 2^{10}$

Even with $2^{10}$ sample points, the sample means of the inner products $\langle\chi_i, \chi_j\rangle$ approximate very well to their expected values. This comparison shows that the orthogonality relations of irreducible characters is much more suitable for the study of Sato-Tate groups than using moment sequences.

[†] We use the fact that all the irreducible characters of $\mathrm{USp}(2g)$ are real.

EXAMPLE 3. We consider the family of non-hyperelliptic genus 3 curves $C$ with an involution. This family is studied in the first part of the author's thesis [13]. Generically, such curve $C$ admits an affine form

$$C : y^4 + g(x)\, y^2 + h(x) = 0$$

with $\deg_x(g) = 2$ and $\deg_x(g) = 4$. The involution gives a degree two map $C \to E$ to an elliptic curve $E$, and thus an isogenous decomposition $0 \to A \to \mathrm{Jac}(C) \to E \to 0$. The image of the Frobenius $\mathrm{Frob}_C(\mathfrak{p})$ on $\mathrm{Jac}(C)$ under

$$\mathrm{Aut}(V_\ell(\mathrm{Jac}(C))) \xrightarrow{\sim} \mathrm{Aut}(V_\ell(E \times A)) \xrightarrow{\sim} \mathrm{Aut}(V_\ell(E) \times V_\ell(A))$$

is $(\mathrm{Frob}_E(\mathfrak{p}), \mathrm{Frob}_A(\mathfrak{p})) \in \mathrm{Aut}(V_\ell(E)) \times \mathrm{Aut}(V_\ell(A))$. We first study the Frobenius distribution of $\mathrm{Frob}_E(\mathfrak{p})$ and $\mathrm{Frob}_A(\mathfrak{p})$ over the family respectively. We compute the data for $p \leq 47$ and over a set of ($\sim 47000$) curves in this family. We use $G = H = \mathrm{SU}(2)$ for $\mathrm{Frob}_E(\mathfrak{p})$ and $G = H = \mathrm{USp}(4)$ for $\mathrm{Frob}_A(\mathfrak{p})$, and the results in Table 5 suggest that both distributions are the generic cases.

| | | | | | |
|---|---|---|---|---|---|
| 1.00 | 0.07 | −0.01 | 0.00 | 0.00 | 0.00 |
| 0.07 | 0.99 | 0.07 | −0.01 | 0.00 | 0.00 |
| −0.01 | 0.07 | 0.99 | 0.07 | −0.01 | 0.00 |
| 0.00 | −0.01 | 0.07 | 0.99 | 0.07 | −0.01 |
| 0.00 | 0.00 | −0.01 | 0.07 | 0.99 | 0.06 |
| 0.00 | 0.00 | 0.00 | −0.01 | 0.06 | 0.92 |

Using $H = \mathrm{SU}(2)$ for $\mathrm{Frob}_E(\mathfrak{p})$

| | | | | | |
|---|---|---|---|---|---|
| 1.00 | 0.00 | −0.06 | 0.07 | 0.00 | 0.02 |
| 0.00 | 1.01 | 0.00 | 0.00 | 0.01 | 0.00 |
| −0.06 | 0.00 | 1.09 | −0.01 | 0.00 | −0.09 |
| 0.07 | 0.00 | −0.01 | 1.02 | 0.00 | 0.07 |
| 0.00 | 0.01 | 0.00 | 0.00 | 1.06 | 0.00 |
| 0.02 | 0.00 | −0.09 | 0.07 | 0.00 | 1.08 |

Using $H = \mathrm{USp}(4)$ for $\mathrm{Frob}_A(\mathfrak{p})$

Table 5

Now we study the Frobenius distribution of $\mathrm{Frob}_C(\mathfrak{p})$ over the family. We guess its Sato-Tate group is $H = \mathrm{SU}(2) \times \mathrm{USp}(4)$. Using the irreducible characters of $G = \mathrm{USp}(6)$, we obtain

| | | | | | |
|---|---|---|---|---|---|
| 1.00 | 0.07 | 0.94 | −0.27 | 0.07 | −0.16 |
| 0.07 | 2.01 | −0.35 | 0.95 | −0.06 | 2.04 |
| 0.94 | −0.35 | 3.00 | −0.31 | 1.06 | −1.16 |
| −0.27 | 0.95 | −0.31 | 2.13 | −0.70 | 2.05 |
| 0.07 | −0.06 | 1.06 | −0.70 | 3.07 | −1.16 |
| −0.16 | 2.04 | −1.16 | 2.05 | −1.16 | 6.24 |

Using $H = \mathrm{SU}(2) \times \mathrm{USp}(4) \subset G = \mathrm{USp}(6)$ for $\mathrm{Frob}_C(\mathfrak{p})$

| | | | | | |
|---|---|---|---|---|---|
| 1 | 0 | 1 | 0 | 0 | 0 |
| 0 | 2 | 0 | 1 | 0 | 2 |
| 1 | 0 | 3 | 0 | 1 | −1 |
| 0 | 1 | 0 | 2 | −1 | 2 |
| 0 | 0 | 1 | −1 | 3 | −1 |
| 0 | 2 | −1 | 2 | −1 | 6 |

Rounded values

| | | | | | |
|---|---|---|---|---|---|
| 1 | 0 | 1 | 0 | 0 | 0 |
| 0 | 2 | 0 | 1 | 0 | 2 |
| 1 | 0 | 3 | 0 | 1 | 0 |
| 0 | 1 | 0 | 2 | 0 | 2 |
| 0 | 0 | 1 | 0 | 3 | 0 |
| 0 | 2 | 0 | 2 | 0 | 6 |

Expected values

This suggests that $\mathrm{SU}(2) \times \mathrm{USp}(4)$ should be the Sato-Tate group, despite that we obtain some entries with value $-1$ in the rounded values, which are caused by the small number of primes used to produce the sample.

It is clear that the Sato-Tate group is contained in $G = \mathrm{SU}(2) \times \mathrm{USp}(4)$, which is the smallest group that we know (for free) containing the Sato-Tate group. The conjectural Sato-Tate group $H$ is $G$ itself. Instead of using the irreducible characters of $\mathrm{USp}(6)$, we use the irreducible characters of $G$, which are products of the irreducible characters of $\mathrm{SU}(2)$ and $\mathrm{USp}(4)$ respectively. We take the first four irreducible characters from each factor to form a set of 16 irreducible characters of $G$. We expect to obtain orthonormal relations, and Table 7 supports our conjecture with very good approximations.

| | | | | | | | | | | | | | | | |
|---|---|---|---|---|---|---|---|---|---|---|---|---|---|---|---|
| 1 | 0.1 | 0 | 0 | 0 | 0 | 0 | 0 | −0.1 | −0.3 | 0 | 0 | 0.1 | 0 | 0 | 0 |
| 0.1 | 1 | 0.1 | 0 | 0 | 0 | 0 | 0 | −0.3 | −0.1 | −0.3 | 0 | 0 | 0.1 | 0 | 0 |
| 0 | 0.1 | 1 | 0.1 | 0 | 0 | 0 | 0 | 0 | −0.3 | −0.1 | −0.3 | 0 | 0 | 0.1 | 0 |
| 0 | 0 | 0.1 | 1 | 0 | 0 | 0 | 0 | 0 | 0 | −0.3 | −0.1 | 0 | 0 | 0 | 0.1 |
| 0 | 0 | 0 | 0 | 1 | −0.2 | 0 | 0 | 0 | 0 | 0 | 0 | 0 | 0 | 0 | 0 |
| 0 | 0 | 0 | 0 | −0.2 | 1 | −0.2 | 0 | 0 | 0 | 0 | 0 | 0 | 0 | 0 | 0 |
| 0 | 0 | 0 | 0 | 0 | −0.2 | 1 | −0.2 | 0 | 0 | 0 | 0 | 0 | 0 | 0 | 0 |
| 0 | 0 | 0 | 0 | 0 | 0 | −0.2 | 1 | 0 | 0 | 0 | 0 | 0 | 0 | 0 | 0 |
| −0.1 | −0.3 | 0 | 0 | 0 | 0 | 0 | 0 | 1 | 0.1 | 0 | 0 | 0 | −0.3 | 0 | 0 |
| −0.3 | −0.1 | −0.3 | 0 | 0 | 0 | 0 | 0 | 0.1 | 1 | 0.1 | 0 | −0.3 | 0 | −0.2 | 0 |
| 0 | −0.3 | −0.1 | −0.3 | 0 | 0 | 0 | 0 | 0 | 0.1 | 1 | 0.1 | 0 | −0.2 | 0 | −0.2 |
| 0 | 0 | −0.3 | −0.1 | 0 | 0 | 0 | 0 | 0 | 0 | 0.1 | 1 | 0 | 0 | −0.2 | 0 |
| 0.1 | 0 | 0 | 0 | 0 | 0 | 0 | 0 | 0 | −0.3 | 0 | 0 | 1 | −0.2 | 0 | 0 |
| 0 | 0.1 | 0 | 0 | 0 | 0 | 0 | 0 | −0.3 | 0 | −0.2 | 0 | −0.2 | 1 | −0.1 | 0 |
| 0 | 0 | 0.1 | 0 | 0 | 0 | 0 | 0 | 0 | −0.2 | 0 | −0.2 | 0 | −0.1 | 1 | −0.1 |
| 0 | 0 | 0 | 0.1 | 0 | 0 | 0 | 0 | 0 | 0 | −0.2 | 0 | 0 | 0 | −0.1 | 1 |

Table 7: Using $H = G = \mathrm{SU}(2) \times \mathrm{USp}(4)$ for $\mathrm{Frob}_C(\mathfrak{p})$

Example 4. We study the curve $C : y^2 = x^8 + 1$, which is studied in [3]. We have the quotient maps:

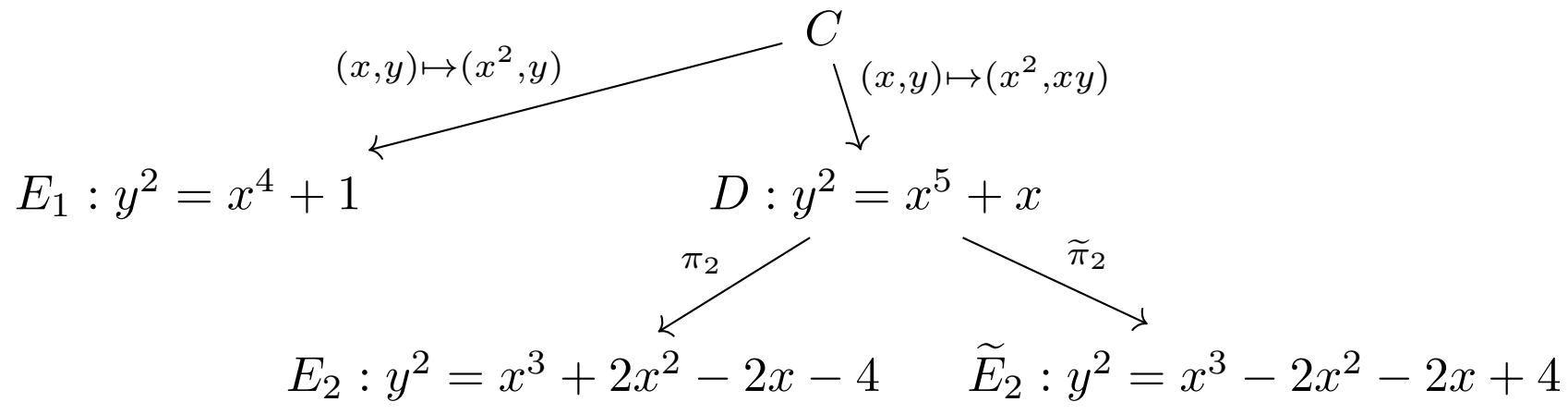

where $\pi_2 : (x, y) \mapsto (x + \frac{1}{x}, y(\frac{1}{x} + \frac{1}{x^2}))$, $\widetilde{\pi}_2 : (x, y) \mapsto (x + \frac{1}{x}, y(\frac{1}{x} - \frac{1}{x^2}))$. The two elliptic curves $E_2$ and $\widetilde{E}_2$ are isogenous, and we have a isogenous decomposition of $\mathrm{Jac}(C) \sim E_1 \times E_2 \times \widetilde{E}_2$. We focus on the identity component, and thus we restrict ourselves to primes $p \equiv 1 \pmod 8$, for which $E_2 \cong \widetilde{E}_2$ over $\mathbb{F}_p$. The identity component of the Sato-Tate group of $C$ is determined by $E_1$ and $E_2$, which both have CM, hence it should be $H = \mathrm{SO}(2)^2$. We work with the irreducible characters of $G = \mathrm{USp}(4)$, and we expect the orthogonality relations are determined by the branching rule from $G$ to $H$. We verify this fact using $2^{12}$ sample points for $E_1 \times E_2$, rather than Frobenius sample points of $C$ for $p \le 2^{40}$ in [3].

| | | | | |
|---|---|---|---|---|
| 1.00 | 0.01 | 0.98 | 1.97 | 0.07 |
| 0.01 | 3.95 | 0.09 | 0.16 | 7.75 |
| 0.98 | 0.09 | 4.92 | 5.80 | 0.30 |
| 1.97 | 0.16 | 5.80 | 11.60 | 0.62 |
| 0.07 | 7.75 | 0.30 | 0.62 | 22.9 |

Using $H = \mathrm{SO}(2)^2 \subset G = \mathrm{USp}(4)$

| | | | | |
|---|---|---|---|---|
| 1 | 0 | 1 | 2 | 0 |
| 0 | 4 | 0 | 0 | 8 |
| 1 | 0 | 5 | 6 | 0 |
| 2 | 0 | 6 | 12 | 0 |
| 0 | 8 | 0 | 0 | 24 |

Expected values

Example 5. We study how well the sample means approximate the expected values in the sample size $n$, which is measured by the function Err defined after Equation (5.1). We choose the elliptic curve $E : y^2 = x^3 + x + 1$. Its Sato-Tate group is $H = \mathrm{SU}(2)$, and we use the set $I$ of its first nine dominant weights. We compute $2^{26}$ Frobenius and plot the function $\delta(k) := \mathrm{Err}(H, S_n, I)$, where $S_n$ is the set of the first $n$ primes of good reduction of $E$, for $n = 2^{10}k$, $k = 1, 2, \ldots, 2^{16}$. The following table shows the pictures of $\delta(k)$ and $\bar{\delta}(k) := \sqrt{\frac{\sum_{i=1}^{k} \delta(i)^2}{k}}$:

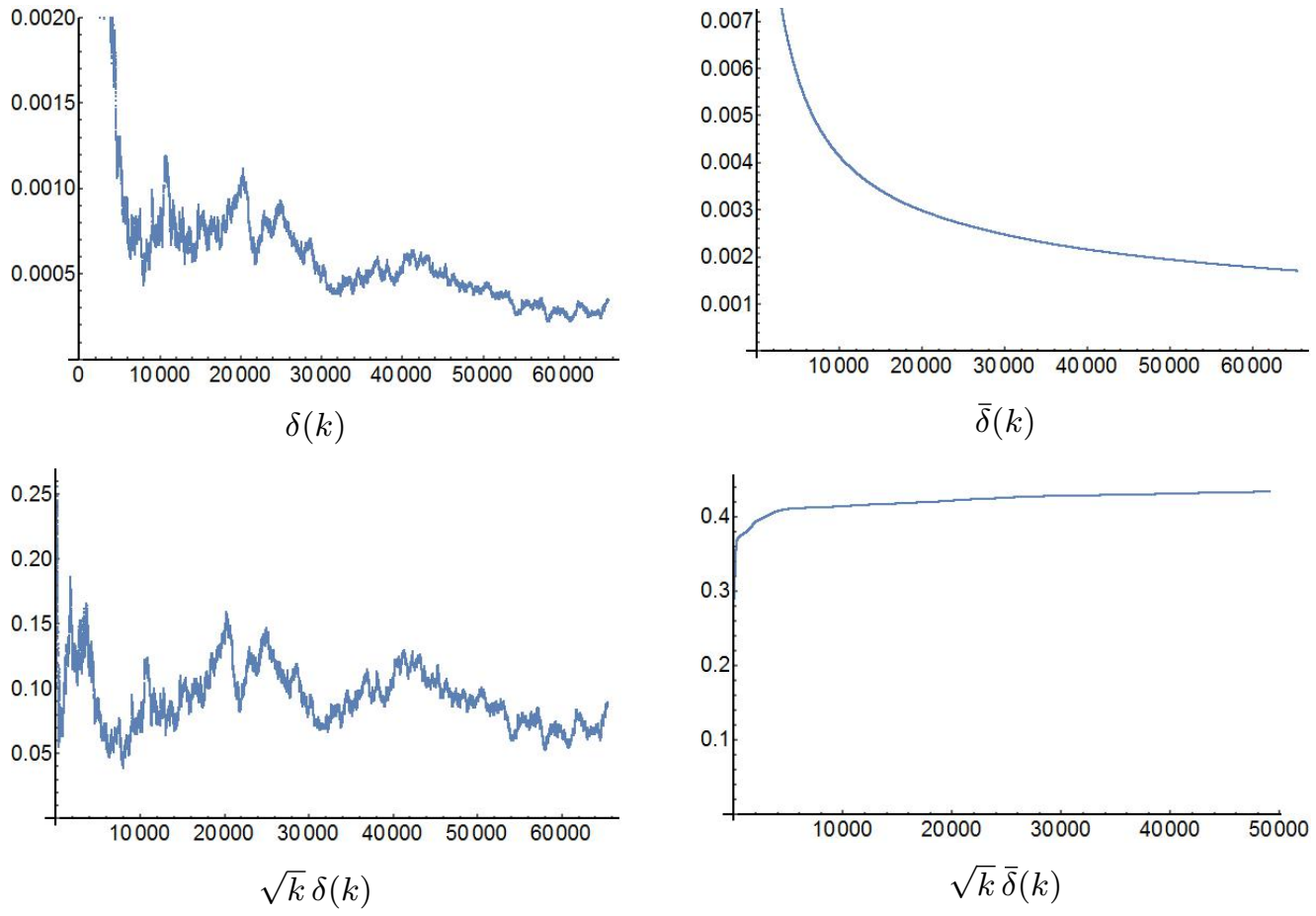

The pictures on the left side have oscillation, but the picture of $\sqrt{k}\,\bar{\delta}(k)$ suggests that it converges to a constant $c$. After a change to the variable $n$, one guess $\operatorname{Err}(H, S_n, I) \approx \frac{32c}{\sqrt{n}}$.

## 6. *Conclusion*

We have developed a systematic way for the computation of the irreducible characters of $\mathrm{USp}(2g)$ in terms of the coefficients $s_i$ of the real Weil polynomial $g_{\mathfrak{p}}$. The main tool is the Brauer-Klimyk formula (Theorem 4.1). We obtain the recursive Algorithm 1 for the computation of the irreducible characters. Although we work with $\mathrm{USp}(2g)$, the algorithm can be modified to compute the irreducible characters of other compact connected Lie groups. In fact, the Brauer-Klimyk formula is already used in Sage to decompose tensor products of two irreducible representations into direct sum of irreducible representations, and it works with a wide collection of classical and exceptional Lie groups (see the Sage documentation [2]). However, using the Brauer-Klimyk formula in the form of Algorithm 1 is new.

The use of orthogonality relations of irreducible characters provides a new perspective to the study of Sato-Tate groups. In Example 2, we show that our new approach requires many fewer sample points to identify the Sato-Tate group $\mathrm{USp}(4)$ in contrast to the approach using moment sequences (in fact $2^{10}$ to $2^{12}$ points provide satisfactory convergence). In Example 3 and Example 4, we demonstrate that it is better to use the character theory of the smallest group we know containing the Sato-Tate group. When we study families of curves with particular structures, like RM curves, this is very useful. This way, the orthogonality relations are small integers. Combine the results in Example 5, we believe that a small size of sample points is enough not only for the generic case, but also for all of the possible connected Sato-Tate groups (or their identity components). Furthermore, in the cases where we don't know the structure of a target curve or family, we can start with the irreducible characters of $\mathrm{USp}(2g)$. It is very likely that we get useful information from them, but without a very good convergence to distinguish the Sato-Tate group from just a few possible candidates. Then we use the character theory for these possible groups to find out the actual one.

We have established the general framework and necessary tools in this article for the study of Sato-Tate groups using the orthogonality relations of irreducible characters. We have seen the heuristic advantages of this new method. We focused on the analysis of the identity component of the Sato-Tate group. A complete strategy requires an analysis of the component group, which is a finite Galois group. A similar approach through character theory of finite group should have as the goal to determine its splitting field, which is usually determined by the geometry and arithmetic of the curve or abelian variety, particularly the splitting field for its automorphism group or endomorphism ring.

Further studies of Sato-Tate groups using this method, in particular for $g = 2$ and $g = 3$, is under way.

*Yih-Dar SHIEH*
*Institut de Mathématiques de Marseille*
*163, avenue de Luminy, Case 907*
*13288 Marseille Cedex 9*
*France*

chiapas@gmail.com

*Author's webpage*